\documentclass[twoside,leqno]{article}
\usepackage[centertags]{amsmath}
\usepackage{amssymb}
\usepackage{amsthm}
\usepackage{hyperref}

\newtheorem{thm}{Theorem}[section]

\newtheorem{prop}[thm]{Proposition}
\newtheorem{defn}[thm]{Definition}

\numberwithin{equation}{section}

\newcommand{\condexp}[1]{\left|}
\newcommand{\Real}{\mathbb R}

\newcommand{\eps}{\epsilon}

\newcommand{\zbs}{\boldsymbol{z}}
\newcommand{\To}{\longrightarrow}
\newcommand{\mbs}{\boldsymbol{m}}

\newcommand{\X}{\boldsymbol{X}}
\newcommand{\Y}{\boldsymbol{Y}}
\newcommand{\A}{\mathcal{A}}
\newcommand{\punt}{\mathbf{.}}
\newcommand{\trasp}{\footnotesize T}

\newcommand{\nubs}{\boldsymbol{\nu}}
\newcommand{\chibs}{\boldsymbol{\chi}}
\newcommand{\gabs}{\boldsymbol{\eta}}
\newcommand{\ibs}{\boldsymbol{i}}
\newcommand{\deltabf}{\boldsymbol{\varsigma}}
\newcommand{\iotabf}{\boldsymbol{\iota}}
\newcommand{\mubf}{\boldsymbol{\mu}}
\newcommand{\vbs}{\boldsymbol{v}}
\newcommand{\jbs}{\boldsymbol{j}}
\newcommand{\varbs}{\boldsymbol{\epsilon}}
\newcommand{\lambs}{\boldsymbol{\lambda}}

\newcommand{\kbs}{\boldsymbol{k}}
\newcommand{\xbs}{\boldsymbol{x}}

\newcommand{\Ybs}{\boldsymbol{Y}}
\newcommand{\ubs}{\boldsymbol{u}}
\newcommand{\boeta}{\boldsymbol{\eta}}

\begin{document}
\title{On a representation of time space-harmonic polynomials via symbolic L\'evy processes.}
\author{Elvira Di Nardo}
\date{}
\maketitle
\vspace{-0.5cm}
\begin{center}
{\sl \small{Dept. Mathematics, Computer Science and Economics, \\
Univ. Basilicata, Viale dell'Ateneo Lucano 10, 85100, Potenza, \\
elvira.dinardo@unibas.it}}
\end{center}
\begin{abstract}
In this paper, we review the theory of time space-harmonic polynomials developed by using a 
symbolic device known in the literature as the classical umbral calculus. The advantage of  this symbolic tool is twofold. First a moment representation is allowed for a wide class 
of polynomial stochastic involving the L\'evy processes in respect to which they are martingales. This representation includes some well-known examples such as Hermite polynomials in connection with Brownian motion. As a consequence, characterizations of many other families of polynomials having the time space-harmonic property can be recovered via the symbolic moment representation. New relations with Kailath-Segall polynomials are  stated.  Secondly the generalization to the multivariable framework is straightforward. Connections with cumulants and Bell polynomials are highlighted both in the univariate case and in the multivariate one. Open problems are addressed at the end of the paper.
\end{abstract}
{\bf keywords:} L\'evy process, time-space harmonic polynomial, Kailath-Segall polynomial, cumulant, umbral calculus. 
%
\section{Introduction}
In mathematical finance, a L\'evy process \cite{sato} is usually employed to model option pricing. 
\begin{defn}
A L\'evy process $X=\{X_t\}_{t \geq 0}$ is a stochastic process satisfying the following properties:
\begin{description}
\item[{\it a)}] $X$ has independent and stationary increments;
\item[{\it b)}] $P[X(0)=0]=1$ on the probability space $(\Omega, {\cal F}, P);$ 
\item[{\it c)}] $X$ is stochastically continuous, i.e. for all $a > 0$ and for all $s \geq 0,$ 
$\lim_{t \rightarrow s} P(|X(t)-X(s)|>a)=0.$ 
\end{description}
\end{defn}
The employment of L\'evy processes in mathematical finance is essentially due to the property of manage continuous 
processes interspersed with jump discontinuities of random size and at random times, well fitting the
main dynamics of a market. In order to include the risk neutrality, a martingale pricing could be
applied to options. But L\'evy processes do not necessarily share the martingale property
unless they are centred. Instead of focusing the attention
on the expectation, a different approach consists in resorting a family of stochastic processes, called {\it polynomial processes} and introduced very recently in \cite{cuchiero}. These processes are built by considering a suitable family of polynomials $\{P(x,t)\}_{t \geq 0}$ and by replacing the indeterminate $x$ with a stochastic process $X_t.$
Introduced in \cite{Seng00} and called {\it time-space harmonic polynomials} (TSH), the polynomials $\{P(x,t)\}_{t \geq 0}$ are such that 
\begin{equation}
E[P(X_t,t) \,\, |\; \mathfrak{F}_{s}] =P(X_s,s), \qquad \hbox{for} \,\, s \leq t
\label{TSH}
\end{equation}
where $\mathfrak{F}_{s}=\sigma\left( X_\tau : \tau \leq s\right)$ 
is the natural filtration associated with $\{X_t\}_{t \geq 0}.$

As done in \cite{neveu} for the discretized version of a L\'evy process, that is a random walk, 
TSH polynomials can be characterized as coefficients of the Taylor expansion
\begin{equation}
\frac{\exp\{z X_t\}}{E[\exp\{z X_t\}]} = \sum_{k \geq 0} R_k(X_t,t)\frac{z^k}{k!}
\label{expmart}
\end{equation}
in some neighborhood of the origin. The left-hand side of (\ref{expmart}) is the so-called 
Wald's exponential martingale \cite{Kuchler}. Wald's exponential martingale is well 
defined only when the process admits moment generating
function $E[\exp\{z X_t\}]$ in a suitable neighborhood of the
origin. Different authors have tried to overcome this gap by
using other tools.  Sengupta \cite{Seng00} uses a
discretization procedure to extend the results proved by Goswami
and Sengupta in \cite{GS95}. Sol\'e and Utzet \cite{Sole} use Ito's formula
showing that TSH polynomials with respect to
L\'evy processes are linked to the exponential complete Bell 
polynomials \cite{Comtet}. Wald's exponential martingale
(\ref{expmart}) has been recently reconsidered also in
\cite{Sengupta08}, but without this giving rise to a closed expression
for these polynomials. 

The employment of the classical umbral calculus turns out to be crucial in dealing with 
(\ref{expmart}). Indeed, the expectation of the polynomial processes
$R_k(X_t,t)$ can be considered without taking into account any question involving the convergence of the right hand side of (\ref{expmart}). 
Indeed the family $\{R_k(x,t)\}_{t \geq 0}$ is linked to the Bell polynomials which are one of the building 
blocks of the symbolic method. The main point here is that any TSH polynomial could be expressed 
as a linear combination  of the family $\{R_k(x,t)\}$ and the symbolic representation of these coefficients 
is particularly suited to be implemented in any symbolic software. The symbolic approach highlights the 
role played by L\'evy processes with regard to which the property (\ref{TSH}) holds and makes clear the
dependence of this representation on their cumulants. 

The paper is organized as follows. Section 2 is provided for readers unaware of the classical umbral calculus. We 
have chosen to recall terminology, notation and the basic definitions strictly necessary to deal with the object of 
this paper. We skip any proof. The reader interested in is referred to \cite{Dinsen, Dinardoeurop}. The theory of 
TSH polynomials is resumed in Section 3 together with the symbolic representation of L\'evy 
processes closely related to their infinite divisible property. Umbral expressions of many classical families of 
polynomials as TSH polynomials with respect to suitable L\'evy processes are outlined. The 
generalization to the multivariable framework is given in Section 4. This setting allows us to deal with multivariate 
Hermite, Euler and Bernoulli polynomials as well as with the class of multivariate L\'evy-Sheffer systems introduced in
\cite{DiNardoOliva3}. Open problems are addressed at the end of the paper. 
\section{The classical umbral calculus.}
Let ${\mathbb R}[x]$ be the ring of polynomials  with real coefficients\footnote{The ring ${\mathbb R}[x]$ may be replaced by any ring in whatever number of indeterminates, as for example
${\mathbb R}[x, y, \ldots].$} in the indeterminate $x.$ The classical umbral calculus is a syntax consisting in an alphabet $\A=\{\alpha,\beta,\gamma, \ldots\}$
of elements, called \emph{umbrae}, and a linear functional $E\,:\,{\mathbb R}[x][\A]\To {\mathbb R}[x]$,
called \emph{evaluation}, such that $E[1]=1$ and
$$E[x^n \, \alpha^i \, \beta^j \, \cdots \, \gamma^k ]=x^n \, E[\alpha^i] \, E[\beta^j] \, \cdots \,  E[\gamma^k] 
\quad \hbox{(uncorrelation property)}$$
where $\alpha, \beta, \ldots, \gamma$ are distinct umbrae and $n,i,j, \ldots,k$ are nonnegative integers.

A sequence $\{a_i\}_{i \geq 0} \in {\mathbb R}[x],$ with $a_0=1,$  is \emph{umbrally represented}\footnote{When no misunderstanding occurs, we use the notation $\{a_i\}$ instead of $\{a_i\}_{i \geq 0}$} by an umbra
$\alpha$ if $E[\alpha^i]=a_i,\;$ for all nonnegative integers $i.$  Then $a_i$ is called the $i$-th
{\it moment} of $\alpha$. An umbra is \textit{scalar} if its moments are elements of ${\mathbb R}$
while it is \textit{polynomial} if its moments are polynomials of ${\mathbb R}[x].$ Special scalar 
umbrae are given in Table 1.
\begin{table}[h]
\centering
\renewcommand{\arraystretch}{1.1}
\begin{tabular}{|c|l|} \hline
\emph{Umbrae} & \emph{Moments} \\  \hline
Augmentation $\epsilon$ & $E[\epsilon^i]=\delta_{0,i},$ with $\delta_{i,j}=1$ if $i=j,$ otherwise $\delta_{i,j}=0.$ \\ \hline
Unity $u$    & $E[u^i]=1$  \\ \hline
Boolean unity $\bar{u}$ & $E[\bar{u}^i] = i!$  \\ \hline
Singleton $\chi$ & $E[\chi]=1$ and $E[\chi^i]=0$, for all $i > 1$ \\ \hline
Bell $\beta$ & $E[\beta^i]=B_i,$ with $B_i$ the $i$-th Bell number \\ \hline
Bernoulli $\iota$ & $E[\iota^i]={\mathfrak B}_i,$  with ${\mathfrak B}_i$ the $i$-th Bernoulli number \\  \hline
Euler $\varepsilon$ & $E[\varepsilon^i]={\mathfrak E}_i,$  with ${\mathfrak E}_i$ the $i$-th Euler number \\ \hline
\end{tabular}
\caption{Special scalar umbrae. The equalities on the right column refer to all nonnegative integer $i,$ unless otherwise specified.}
\label{table1}
\end{table}
The core of this moment symbolic calculus consists in defining the {\it dot-product} of two umbrae,
whose construction is shortly recalled in the following. 

First let us underline that in the alphabet $\A$ two (or more) distinct umbrae may represent the same sequence of moments. More formally, two umbrae $\alpha$ and $\gamma$ are said to be {\it
similar} when $E[\alpha^n]=E[\gamma^n]$ for all nonnegative integers $n,$ in symbols
$\alpha \equiv \gamma.$ Therefore, given a sequence $\{a_n\},$ there are infinitely many distinct, and thus similar umbrae, representing the sequence. 

Denote $\alpha'+\alpha''+\cdots+\alpha'''$ by the symbol $n \punt \alpha,$ where $\{\alpha',\alpha'',\ldots,\alpha'''\}$ is a set of $n$ uncorrelated umbrae similar to $\alpha.$ 
The symbol $n \punt \alpha$ is an example of \emph{auxiliary umbra}. In a \emph{saturated} umbral calculus, the auxiliary umbrae are managed as they were elements of $\A$ \cite{SIAM}.  The umbra $n \punt \alpha$ is called the 
{\it dot-product} of the integer $n$ and the umbra $\alpha$ with moments 
\cite{Dinardoeurop}:
\begin{equation}
q_i(n)=E[(n \punt \alpha)^i]=\sum_{k=1}^i (n)_k B_{i,k}(a_1, a_2, \ldots, a_{i-k+1}),
\label{(1)}
\end{equation}
where $(n)_k$ is the lower factorial and $B_{i,k}$ are the exponential partial Bell polynomials
\cite{Comtet}.

In (\ref{(1)}), the polynomial $q_i(n)$ is of degree $i$ in $n.$ If the integer $n$ is replaced by
$t \in {\mathbb R},$ in (\ref{(1)}) we have $q_i(t) = \sum_{k=1}^i (t)_k
B_{i,k}(a_1, a_2, \ldots, a_{i-k+1}).$ Denote by $t \punt \alpha$ the auxiliary umbra 
such that $E[(t \punt \alpha)^i] = q_i(t),$ for all nonnegative integers $i$.
The umbra $t \punt \alpha$ is the dot-product of $t$ and $\alpha.$ A kind of distributive property holds:
\begin{equation}
(t + s) \punt \alpha \equiv t \punt \alpha + s \punt \alpha^{\prime}, \quad s,t \in {\mathbb R}
\label{(distributive)}
\end{equation}
where $\alpha^{\prime} \equiv \alpha.$ In particular if in (\ref{(1)}) the integer $n$ is replaced by $-t,$ the auxiliary umbra $-t \punt \alpha$ is such that
\begin{equation}
-t \punt \alpha + t \punt \alpha^{\prime} \equiv \eps,
\label{(inverse)}
\end{equation}
where $\alpha^{\prime} \equiv \alpha.$ Due to equivalence (\ref{(inverse)}), the umbra $-t \punt \alpha$ is the {\it inverse}\footnote{Since $-t \punt \alpha$ and $t \punt \alpha$ are two distinct symbols, they are considered uncorrelated, therefore $-t \punt \alpha + t \punt \alpha^{\prime} \equiv 
-t \punt \alpha + t \punt \alpha \equiv \eps.$ When no confusion occurs, we will use this last similarity instead of
(\ref{(inverse)}).} umbra of $t \punt \alpha.$ 

Let us consider again the polynomial $q_i(t)$ and suppose to replace $t$ by an umbra $\gamma.$ The
polynomial $q_i(\gamma) \in {\mathbb R}[x][\A]$ is an {\it umbral polynomial}  with support 
\footnote{The support $\hbox{\rm supp} \, (p)$ of an umbral polynomial
$p \in {\mathbb R}[x][\A]$ is the set of all umbrae occurring in it.} $\hbox{\rm supp} \,
(q_i(\gamma))=\{\gamma\}.$ The {\it dot-product of $\gamma$ 
and $\alpha$} is the auxiliary umbra $\gamma \punt \alpha$ such that $E[(\gamma \punt \alpha)^i]=E[q_i(\gamma)]$ for all
nonnegative integers $i.$ Two umbral polynomials $p$ and $q$ are said to be \emph{umbrally equivalent} 
if $E[p]=E[q],$ in symbols $p \simeq q.$ Therefore equations (\ref{(1)}), with $n$ replaced by an umbra $\gamma,$ 
can be written as the equivalences
\begin{equation}
q_i(\gamma) \simeq (\gamma \punt \alpha)^i \simeq \sum_{k=1}^i (\gamma)_k B_{i,k}(a_1, a_2, \ldots, a_{i-k+1}).
\label{(dotprodumbrae)}
\end{equation}
Special dot-product umbrae are the  $\alpha$-cumulant umbra $\chi \punt \alpha$ and $\alpha$-partition
umbra $\beta \punt \alpha,$ that we will use later on. In particular any umbra is a partition umbra \cite{Dinardoeurop}. This property means that 
if $\{a_i\}$ is a sequence umbrally represented by an umbra $\alpha,$
then there exists a sequence $\{h_i\}$ umbrally represented by an umbra
$\kappa_{{\scriptscriptstyle \alpha}},$ such that $\alpha \equiv \beta \punt
\kappa_{{\scriptscriptstyle \alpha}}.$ The umbra $\kappa_{{\scriptscriptstyle \alpha}}$ is similar to the $\alpha$-cumulant umbra, that is $\kappa_{{\scriptscriptstyle \alpha}} \equiv \chi \punt \alpha,$
and its moments share the well-known properties of cumulants. \footnote{For cumulants $\{C_i(Y)\}$ of a random variable $Y,$ the following properties hold for all nonnegative integers $i:$ (Homogeneity) $C_i(a Y) = a^i C_i (Y)$ for $a \in \Real,$ 
(Semi-invariance) $C_1(Y+a)=a + C_1(Y), C_i(Y+a) = C_i(Y)$ for $i \geq 2,$ (Additivity) $C_i(Y_1+Y_2) = C_i(Y_1) + C_i(Y_2),$ if $Y_1$ and $Y_2$ are independent random variables.}

Dot-products can be nested. For example, moments of $(\alpha \punt \varsigma) \punt \gamma$
can be recursively computed by applying two times formula (\ref{(dotprodumbrae)}). Parenthesis can be avoided
since $(\alpha \punt \varsigma) \punt \gamma \equiv \alpha \punt (\varsigma \punt \gamma).$
In particular $\alpha \punt \beta \punt \gamma,$ with $\beta$ the Bell umbra, is the so-called {\sl composition umbra}, with moments 
\begin{equation}
E[(\alpha \punt \beta \punt \gamma)^i]=\sum_{k=1}^i a_k B_{i,k}(g_1, g_2, \ldots, g_{i-k+1}),
\label{(2)}
\end{equation}
where $\{a_i\}$ are moments of $\alpha$ and $\{g_i\}$ are moments of $\gamma.$ When the umbra $\alpha$ is replaced by $t \in \Real,$ then equation (\ref{(2)}) gives 
the $i$-th moment of a compound Poisson random variable (r.v.) of parameter $t:$
$$E[(t \punt \beta \punt \gamma)^i]=\sum_{k=1}^i t^k B_{i,k}(g_1, g_2, \ldots, g_{i-k+1}).$$
There are more auxiliary umbrae that will employed in the following. 
For example, if $E[\alpha] \ne 0,$ the compositional inverse $\alpha^{\scriptscriptstyle <-1>}$ of an umbra $\alpha$ is such that $\alpha \punt \beta \punt \alpha^{\scriptscriptstyle <-1>} \equiv
\alpha^{\scriptscriptstyle <-1>} \punt \beta \punt \alpha \equiv \chi.$ The derivative of an umbra $\alpha$ is the umbra $\alpha_{\scriptsize D}$ whose moments 
are $E[\alpha_{\scriptsize D}^i] =  i \, a_{i-1}$ for all nonnegative integers $i \geq 1.$ 
The disjoint sum $\alpha \dot{+} \gamma$ of $\alpha$ and $\gamma$  represents the sequence $\{a_i + g_i\}.$ Its main property involves the Bell umbra:
\begin{equation}
\beta \punt (\alpha \dot{+} \gamma) \equiv \beta \punt \alpha + \beta \punt \gamma.
\label{(disj)}
\end{equation}
\subsection{Symbolic L\'evy processes.}
%
The family of auxiliary umbrae $\{t \punt \alpha\}_{t \in I},$ with $I \subset  {\mathbb R}^{+},$
is the umbral counterpart of a stochastic process $\{X_t\}_{t \in I}$ having all moments
and such that $E[X_t^i]=E[(t \punt \alpha)^i]$ for all nonnegative integers $i.$ This symbolic 
representation parallels the well-known infinite divisible property of a L\'evy process,
summarized by the following equality in distribution
\begin{equation}
X_t \stackrel{d}{=} \underbrace{\Delta X_{t/n} + \cdots + \Delta X_{t/n}}_n
\label{(idp)}
\end{equation}
with $\Delta X_{t/n}$ a r.v. corresponding to the increment of the process over
an interval of amplitude $t/n.$ The $n$-fold convolution (\ref{(idp)}) is usually expressed by the product of $n$ times a characteristic function $E[e^{{\mathfrak i} z X_t}] = E[e^{{\mathfrak i} z \Delta X_{t/n}}]^n$
with ${\mathfrak i}$ the imaginary unit. More generally one has 
\begin{equation}
E[e^{{\mathfrak i} z X_t}] = E[e^{{\mathfrak i} z X_1}]^t.
\label{(chfunc)}
\end{equation}
Equation (\ref{(chfunc)}) allows us to show that the auxiliary umbra $t \punt \alpha$ is the symbolic version
of $X_t.$ To this aim we recall that the formal power series 
\begin{equation}
f(\alpha, z) = 1 + \sum_{i \geq 1} a_i \frac{z^i}{i!}
\label{(fps)}
\end{equation}
is the generating function of an umbra $\alpha,$ umbrally representing the sequence $\{a_i\}.$ 
Table 2 shows generating functions for some special auxiliary umbrae introduced in the previous section.
\begin{table}[ht]
\centering
\renewcommand{\arraystretch}{1.1}
\begin{tabular}{|c|l|} \hline
\emph{Umbrae} & \emph{Generating functions} \\ \hline
Augmentation $\epsilon$ & $f(\epsilon,z)=1$ \\ \hline
Unity $u$    & $f(u,z)=e^z$ \\ \hline
Boolean unity $\bar{u}$ & $f(\bar{u},z)= \frac{1}{1-z}$ \\ \hline
Singleton $\chi$ & $f(\chi,z)= 1 + z$ \\ \hline
Bell $\beta$ & $f(\beta,z)=\exp[e^z-1]$  \\ \hline
Bernoulli $\iota$ & $f(\iota,z) = z/(e^z-1)$ \\ \hline
Euler $\eta$ &  $f(\eta,z) = 2 \, e^z / [e^z + 1]$ \\ \hline
dot-product $n \punt \alpha$ & $f(n \punt \alpha,z)= f(\alpha,z)^n$  \\ \hline
dot-product $t \punt \alpha$ & $f(t \punt \alpha,z)= f(\alpha,z)^t$  \\ \hline
dot-product $\gamma \punt \alpha$ & $f(\gamma \punt \alpha,z)= f(\gamma,\log f(\alpha,z))$  \\ \hline
$\alpha$-cumulant  $\chi \punt \alpha$ & $f(\chi \punt \alpha, z) = 1 + \log[f(\alpha,z)]$ \\ \hline
$\alpha$-partition  $\beta \punt \alpha$ & $f(\beta \punt \alpha,z) = \exp[f(\alpha,z) - 1]$ \\ \hline
composition  $\alpha \punt \beta \punt \gamma$ & $f(\alpha \punt \beta \punt \gamma,z) = f[\alpha,f(\gamma,z) - 1]$ \\ \hline
$\alpha$-partition  $t \punt \beta \punt \gamma$ & $f(t \punt \beta \punt \gamma,z) = \exp[t(f(\gamma,z) - 1)]$ \\ \hline
derivative  $\alpha_{\scriptscriptstyle D}$ & $f(\alpha_{\scriptscriptstyle D},z) = 1 + z f(\alpha, z)$ \\ \hline
\end{tabular}
\caption{Generating functions for some special auxiliary umbrae.}
\label{table2}
\end{table}
\par
As for infinitely divisible stochastic processes (\ref{(idp)}), one has $f(t \punt \alpha, z)=f(\alpha,z)^t.$ 
It is well-known that the class of infinitely divisible distributions coincides with the class of limit distributions of
compound Poisson distributions \cite{Feller}. By the symbolic method, any L\'evy process is 
of compound Poisson type \cite{elviraannali}. This result is a direct consequence of the L\'evy-Khintchine 
formula \cite{sato} involving the moment generating function of a L\'evy process. Indeed, if $\phi(z,t)$ denotes the moment generating function of $X_t$ and $\phi(z)$ denotes the moment generating function of $X_1$ then $\phi(z,t)=\phi(z)^t$ from (\ref{(chfunc)}). From the L\'evy-Khintchine formula $\phi(z)=\exp[g(z)],$  with
\begin{equation}
g(z)= z m + \frac{1}{2} s^2 z^2 + \int_{\Real} \left(e^{zx} - 1 - 
z \, x \, {\bf 1}_{ \{ |x| \leq 1\} } \right) {\rm d}(\nu(x)).
\label{(LK)}
\end{equation}
The term $(m,s^2,\nu)$ is the L\'evy triplet and $\nu$ is the L\'evy measure. The function $\phi(z,t)$ shares  
the same exponential form of the moment generating function $f(t \punt \beta \punt \gamma,z)$ of a compound Poisson process, see Table 2. If $\nu$
admits all moments and if $c_0 = m+\int_{\{|x| \geq 1\}} x \, {\rm d}(\nu(x)),$ then the function 
$g(z)$ given in (\ref{(LK)}) has the form 
\begin{equation}
g(z)= c_0 z + \frac{1}{2} s^2 z^2 + \int_{\Real} \left(e^{zx} - 1 - z \, x \right) {\rm d}(\nu(x)).
\label{(LK1)}
\end{equation}
Thanks to (\ref{(LK1)}), the symbolic representation $t \punt \beta \punt \gamma$ of a L\'evy process is such that the umbra $\gamma$ can be further decomposed in a suitable disjoint sum of umbrae.
Indeed, assume 
\begin{quote} 
{\it i)} $\varsigma$ an umbra with generating function $f(\varsigma,z)=1 + z^2/2,$ \\
{\it ii)} $\, \eta$ an umbra with generating function $f(\eta,z) = 
\int_{\Real} \left(e^{zx} - 1 - z \, x \right) {\rm d}(\nu(x)).$
\end{quote}
Then a L\'evy process is umbrally represented by the family 
\begin{equation}
\{t \punt \beta \punt (c_0 \chi \dot{+} 
s \varsigma \dot{+} \eta)\}_{t \geq 0} \quad \hbox{or} \quad \{t \punt \beta \punt (c_0 \chi \dot{+} 
s \varsigma ) + t \punt \beta \punt \eta\}_{t \geq 0},
\label{(symbLev)}
\end{equation}
due to (\ref{(disj)}). Symbolic representation (\ref{(symbLev)}) is in agreement with It\^{o} representation 
$X_t = W_t + M_t + c_0 t$ of a L\'evy process with $W_t + c_0 t$  a Wiener process and $M_t$ a compensated sum of jumps of a Poisson process involving the L\'evy measure. Indeed the Gaussian component is represented by the symbol $t \punt \beta \punt (c_0 \chi \dot{+} s \varsigma)$ as stated in \cite{dinardooliva2009}, with $c_0$ corresponding to the mean and $s^2$ corresponding to the variance. The Poisson component is represented by $t \punt \beta \punt \eta,$ that is $t \punt \beta \punt \eta$ is the umbral counterpart of a random sum $S_N = Y_1 + \cdots + Y_N,$ with $\{Y_i\}$ independent and identically distributed r.v.'s  corresponding to $\eta,$ associated to the L\'evy measure, and $N$ a Poisson r.v. of parameter $t.$ The representation $\{t \punt \beta \punt (c_0 \chi \dot{+} s \varsigma \dot{+} \eta)\}_{t \geq 0}$ shows that the L\'evy process is itself a compound Poisson process with $\{Y_i\}$ corresponding to the disjoint sum $(c_0 \chi \dot{+} s \varsigma \dot{+} \eta).$ 

More insights may be added on the role played by the umbra $c_0 \chi \dot{+} s \varsigma \dot{+} \eta$. Indeed the moment generating function of a L\'evy process can be written as $\phi(z,t)= \exp[t \, g(z)]$ with $g(z)=\log \, \phi(z).$ So the function $g(z)$ in (\ref{(LK1)}) is the cumulant generating function of $X_1$ and $\gamma \equiv c_0 \chi \dot{+} s \varsigma \dot{+} \eta$ is the symbolic representation of a r.v. whose moments are cumulants of $X_1.$ Therefore, in the symbolic representation $t \punt \alpha$ of a L\'evy process, introduced at the beginning of this section, the umbra $\alpha$ is the partition umbra of the cumulant umbra $\gamma \equiv c_0 \chi \dot{+} s \varsigma \dot{+} \eta$ that is $\alpha \equiv \beta \punt \gamma.$ 

This remark suggests the way to construct the boolean and the free version of a L\'evy process by using the boolean and the free cumulant umbra \cite{DinardoPetrulloSenato}. 
\begin{quote} 
{\bf Boolean L\'evy process.} Let $M(z)$ be the ordinary generating function of a r.v. $X,$ that is $M(z)= 1 + \sum_{i \geq 1}  a_i z^i,$ where $a_i = E[X^i]$. The 
boolean cumulants of $X$ are the coefficients $b_i$ of the power series 
$B(z)= \sum_{i \geq 1} b_i z^i$ such that $M(z) = 1/[1-B(z)].$ Denote by $\bar{\alpha}$ the umbra such that
$E[\bar{\alpha}^i] = i! a_i$ for all nonnegative integers $i.$ Then the umbra 
$\varphi_{\scriptscriptstyle\alpha}$ such that $\bar{\alpha} \equiv \bar{u} \punt \beta \punt \varphi_{\scriptscriptstyle\alpha}$ represents the sequence $\{b_i\}$ and is the $\alpha$-boolean 
cumulant umbra. Therefore the symbolic representation of a boolean L\'evy process is 
$t \punt \bar{u} \punt \beta \punt \varphi_{\scriptscriptstyle\alpha}.$  \par
{\bf Free L\'evy process.} The noncrossing (or free) cumulants of $X$ are the
coefficients $r_i$ of the ordinary power series $R(z) = 1 +
\sum_{i \geq 1} r_i z^i$ such that $M(z)=R[z M(z)].$ If $\bar{\alpha}$ 
is the umbra with generating function $M(z),$ then 
the $\bar{\alpha} \,$-free cumulant umbra $\mathfrak{K}_{\scriptscriptstyle \bar{\alpha}}$ represents 
the sequence $\{i! r_i\}.$ Assuming $\bar{\alpha}$ the umbral counterpart of the increment of a L\'evy process over
the interval $[0,1],$ then the symbolic representation of a free L\'evy process is 
$t \punt {\bar{\mathfrak{K}}}_{\scriptscriptstyle \alpha} \punt
\beta \punt {(-1 \punt {\bar{\mathfrak {K}}}_{\scriptscriptstyle
\alpha})_{\scriptscriptstyle D}^{\scriptscriptstyle <-1>}}.$ 
\end{quote}
Some more remarks on the parameters $c_0$ and $s$ may be added. The L\'evy process in (\ref{(symbLev)})  is a martingale if and only if $c_0 = 0,$ see Theorem 5.2.1 in \cite{Applebaum}. When this happens, $E[X_t]=0$ for all $t \geq 0$ and the L\'evy process is said to be centered. Since the parameter $c_0$ allows
the contribution of the singleton umbra $\chi$ in (\ref{(symbLev)}), such an umbra plays a central role in the martingale property of a L\'evy process. If $c_0=0,$ no contribution is given by $\chi$ which indeed does not admit a probabilistic counterpart.

If $s=0,$ the corresponding L\'evy process is a {\it subordinator}, with almost sure non-decreasing paths. The subordinator processes are usually employed to scale the 
time of a L\'evy process. This device is useful to widen or to close the jumps
of the paths in market dynamics. Denote by $T_t$ the subordinator process of $X_t$ chosen independent of $X_t.$ The process $X_{T_t}$ is of L\'evy type too. The symbolic representation of $T_t$ is $t \punt \beta \punt (c_0 \chi \dot{+} \eta^{\prime})$ so that $t \punt \beta \punt (c_0 \chi \dot{+} \eta^{\prime}) \punt \beta \punt (c_0 \chi \dot{+} s \varsigma + \eta)$ represents $X_{T_t}$ with $\eta$ and $\eta^{\prime}$ similar and uncorrelated umbrae. Despite its nested representation, the following result is immediately recovered: the process $X_{T_t}$ is a compound Poisson process $S_N$ with $Y_i$ a randomized compound Poisson r.v. of random parameter $\eta^{\prime},$ shifted of $c_0$ in its mean.  
%
\section{Time-space harmonic polynomials.}
Set ${\cal X}=\{\alpha\}.$ The conditional evaluation $E(\cdot \, \,  \vline \,\, \alpha)$ with respect to $\alpha$ handles the umbra $\alpha$ as it was an indeterminate \cite{elviraannali}. In particular, $E(\cdot \, \,  \vline \,\, \alpha): \Real[x][{\cal A}] \rightarrow \Real[{\cal X}]$ is such that
$E(1 \,\, \vline \,\,\alpha)=1$ and
$$
E(x^m \alpha^n \gamma^i \xi^j\cdots \, \,  \vline \,\, \alpha)=x^m \alpha^n
E[\gamma^i]E[\xi^j]\cdots
$$
for uncorrelated umbrae $\alpha, \gamma, \xi, \ldots$ and for nonnegative integers $m,n,i,j,\ldots.$
As it happens in probability theory, the conditional evaluation is an element of
${\mathbb R}[x][\A]$ and, if we take the overall evaluation of $E(p\,\, \vline \,\, \alpha),$ this gives $E[p \,],$ with $p \in {\mathbb R}[x][\A],$ that is $E[E(p\,\, 
\vline \,\, \alpha)]=E[p \,].$ Umbral polynomials $p,$ not having $\alpha$
in its support, are such that $E(p  \,\, \vline \,\, \alpha)=E[p \,].$ 

Conditional evaluations with respect to auxiliary umbrae need to be handled carefully. For example, since $(n+1) \punt \alpha \equiv n \punt \alpha + \alpha^{\prime},$ the conditional evaluation with respect to the dot product $n \punt \alpha$ is defined as 
$$E[(n+1) \punt \alpha  \,\, \vline \,\, n \punt \alpha] = n \punt \alpha + E[\alpha^{\prime}],$$ 
and more general, from (\ref{(distributive)}) with $t$ and $s$ replaced by $n$ and $m,$ 
\begin{equation}
E( [(n+m) \punt \alpha]^k  \,\, | \,\, n \punt \alpha) =
E([n \punt \alpha + m \punt \alpha^{\prime}]^k \,\, | \,\, n \punt \alpha) = \sum_{j=0}^k \binom{k}{j} (n \punt \alpha)^j E[(m \punt \alpha^{\prime})^{k-j}],
\label{(condevalint)}
\end{equation}
for all nonnegative integers $n$ and $m.$ Therefore, for $t \geq 0$ the conditional evaluation of $t \punt \alpha$ with respect to the auxiliary umbra $s \punt \alpha,$ with $0 \leq s \leq t,$ is defined according to (\ref{(condevalint)}) such as
$$E[(t \punt \alpha)^k  \,\, | \,\, s \punt \alpha] = \sum_{j=0}^k \binom{k}{j} (s \punt \alpha)^j 
E([(t-s) \punt \alpha^{\prime}]^{k-j}).$$
Equation (\ref{TSH}) traces the way to extend the definition of polynomial processes to umbral polynomials. 
\begin{defn}
Let $\{P(x,t)\} \in {\mathbb R}[x][{\cal A}]$ be a family of polynomials indexed by $t \geq 0.$
$P(x,t)$ is said to be a {\rm TSH polynomial} with respect to the family
of umbral polynomials $\{q(t)\}_{t \geq 0}$ if and only if
$E \left[ P(q(t),t) \, \, \vline \, \, q(s) \right] = P(q(s),s)$ for all $0 \leq s \leq t.$
\end{defn}
The main result of this section is the following theorem \cite{elviraannali}.
\begin{thm}\label{UTSH2}
For all nonnegative integers $k,$ the family of polynomials\footnote{When no confusion occurs, we will use the notation $x - t \punt \alpha$ to denote the polynomial umbra $- t \punt \alpha + x = x + (- t) \punt \alpha.$}
$$
Q_k(x,t)=E[(x - t \punt \alpha)^k] \in {\mathbb R}[x]
$$
is TSH with respect to $\{t \punt \alpha\}_{t
\geq 0}.$
\end{thm}
By expanding $Q_k(x,t)$ via the binomial theorem, one has 
$$Q_k(x,t)= \sum_{j=0}^k \binom{k}{j} x^j E[(-t \punt \alpha)^{k-j}]$$
so that
$$Q_k(t \punt \alpha,t)= \sum_{j=0}^k \binom{k}{j} (t \punt \alpha)^j E[(-t \punt \alpha)^{k-j}].$$
Since $t \punt \alpha$ is the symbolic version of a L\'evy process, the property
$$E[Q_k(t \punt \alpha, t) \,\, | \,\, s \punt \alpha] =  
\sum_{j=0}^k \binom{k}{j} E[(t \punt \alpha)^j|s \punt \alpha] E[(-t \punt \alpha)^{k-j}] = 
Q_k(s \punt \alpha, s)$$
parallels equation (\ref{TSH}). In particular $\{Q_k(x,t)\}$ is a polynomial sequence umbrally represented by the polynomial umbra $x - t \punt \alpha,$ which is indeed the TSH polynomial umbra with respect to $t \punt \alpha.$ Polynomial umbrae of type $x + \alpha$ are Appell umbrae \cite{DNS}. Then $\{Q_k(x,t)\}$ is an Appell sequence and
$$\frac{{\rm d}}{{\rm d} x} \, Q_k(x,t) = k \, Q_{k-1}(x,t), \qquad \hbox{for all integers 
$k \geq 1.$}$$ 
The generating function of the TSH polynomial umbra $x -t \punt \alpha$ is 
\begin{equation}
f(x - t \punt \alpha, z)=\frac{\exp\{xz\}}{f(\alpha,z)^t}=\sum_{k \geq 0} Q_k(x,t) \frac{z^k}{k!}.
\label{(wald)}
\end{equation}
By replacing $x$ with $t \punt \alpha$ in (\ref{(wald)}), Wald's exponential martingale (\ref{expmart}) is recovered. Equality of two formal power series is given in terms of equality of their corresponding coefficients, so that $E[R_k(X_t,t)] = E[Q_k(t \punt \alpha,t)]$ by comparing (\ref{(wald)}) with (\ref{expmart}).  Wald's identity $\sum_{k \geq 0} E[R_k(X_t,t)] z^k/k!=1$ is encoded by the equivalence $t \punt \alpha - t \punt \alpha \equiv \epsilon$ obtained from $x -t \punt \alpha$ when $x$ is replaced by $t \punt \alpha.$

The next proposition gives the way to compute the coefficients of $Q_k(x,t)$ in any symbolic 
software.
\begin{prop}
If $\{a_n\}$ is the sequence umbrally represented by the umbra $\alpha$ and 
$\{Q_k(x,t)\}$ is the sequence of TSH polynomials with respect to $\{t \punt \alpha\}_{t \geq 0},$ 
then 
$$Q_k(x,t) = \sum_{j,i=0}^{k} c^{(k)}_{i,j} \, t^i \, x^j,$$
with
$$c^{(k)}_{i,j} = \binom{k}{j} \sum_{\lambda \vdash k-j} {\rm d}_{\lambda} (-1)^{2l(\lambda)+i} \, s[l(\lambda),i]
\, a_1^{r_1} a_2^{r_2} \cdots $$
where the sum is over all partitions\footnote{Recall that a partition of an integer
$i$ is a sequence $\lambda = (\lambda_1, \lambda_2, \ldots,
\lambda_m),$ where $\lambda_j$ are weakly decreasing positive
integers such that $\sum_{j=1}^{m} \lambda_j = i.$ The integers
$\lambda_j$ are named {\it parts} of $\lambda.$ The {\it length}
of $\lambda$ is the number of its parts and will be indicated by
$l(\lambda).$  A different notation is $\lambda = (1^{r_1},
2^{r_2}, \ldots),$ where $r_j$ is the number of parts of $\lambda$
equal to $j$ and $r_1 + r_2 + \cdots = l(\lambda).$ Note that
$r_j$ is said to be the multiplicity of $j$. We use the classical
notation $\lambda \vdash i$ to denote \lq\lq $\lambda$ is a
partition of $i$\rq\rq.} $\lambda = (1^{r_1}, 2^{r_2},
\ldots) \vdash k-j, s[l(\lambda),i]$ 
denotes a Stirling number of first kind and $d_{\lambda} = i!/(r_1!
r_2! \cdots \, (1!)^{r_1}(2!)^{r_2} \cdots).$
\end{prop}
More properties on the coefficients of $Q_k(x,t)$ are given in \cite{elviraannali}. 
\par \medskip
Any TSH polynomial is a linear combination of $\{Q_k(x,t)\},$ which indeed are a bases of the space of TSH polynomials. The following theorem characterizes the coefficients of any TSH polynomial $P(x,t)$ in terms of the coefficients of $\{Q_k(x,t)\}.$
\begin{thm} \label{comb}
A polynomial $P(x,t)=\sum_{j=0}^k p_j(t) \, x^j$ of degree $k$ for all
$t \geq 0$ is a TSH polynomial with respect to $\{t \punt \alpha\}_{t \geq 0}$ 
if and only if
$$p_j(t) = \sum_{i=j}^k  \binom{i}{j} \, p_i(0) \, E[(-t\punt\alpha)^{i-j}], \quad
\hbox{for $j=0,\ldots,k.$}$$
\end{thm}
\subsection{Cumulants.}
A different symbolic representation of TSH polynomials
$\{Q_k(x,t)\}$ is
$$Q_k(x,t)=E[(x - t \punt \beta \punt \kappa_{{\scriptscriptstyle \alpha}})^k],$$
with $\kappa_{{\scriptscriptstyle \alpha}}$ the $\alpha$-cumulant umbra.
The umbra $-t \punt \beta \punt  \kappa_{{\scriptscriptstyle \alpha}} \equiv t \punt 
\beta \punt (-1 \punt \kappa_{{\scriptscriptstyle \alpha}}) \equiv t \punt 
\beta \punt \kappa_{(-1 \punt {\scriptscriptstyle \alpha})}$ is the symbolic version 
of a L\'evy process with sequence of cumulants of $X_1$ umbrally represented by $\kappa_{(-1 \punt {\scriptscriptstyle \alpha})}.$ Therefore, also the polynomials $Q_k(x,t) = E[(x + t \punt \beta \punt \kappa_{{\scriptscriptstyle \alpha}})^k]$ are TSH with respect to L\'evy processes umbrally represented by $\{t \punt 
\beta \punt \kappa_{(-1 \punt {\scriptscriptstyle \alpha})}\}_{t \geq 0} \equiv \{t \punt (-1 \punt \alpha)\}_{t \geq 0}.$

Moments of polynomial umbrae $t \punt \beta \punt \gamma$ involve the exponential Bell polynomials. When $t$ is set equal to $1,$ then complete exponential Bell polynomials
are recovered. More generally, the moments of $x - t \punt \beta \punt \kappa_{{\scriptscriptstyle \alpha}} \equiv x - t \punt \alpha$ can be expressed by using exponential complete Bell polynomials too since
\begin{equation} 
x - t \punt \alpha \equiv \beta \punt [\chi \punt (x - t \punt \alpha)]  
\equiv \beta \punt \kappa_{{\scriptscriptstyle x - t \punt \alpha}}
\equiv \beta \punt (\kappa_{(x \punt u)} \dot{+} \kappa_{(- t \punt \alpha)})
\equiv \beta \punt \kappa_{(x \punt u)} + \beta \punt \kappa_{(- t \punt \alpha)}  
\label{(TSHbell)}
\end{equation}
where $\kappa_{{\scriptscriptstyle x - t \punt \alpha}}$ is the cumulant umbra of
$x - t \punt \alpha,$ that could be replaced by $\kappa_{(x \punt u)} \dot{+} \kappa_{(- t \punt \alpha)}$ due to the additivity property of cumulants. The last equivalence in 
(\ref{(TSHbell)}) follows from equivalence (\ref{(disj)}). From equivalences (\ref{(TSHbell)}), we have
\begin{equation}
Q_k(x,t)=Y_k(x+h_1, h_2, \ldots, h_k),
\label{(complBell)}
\end{equation}
with $Y_k$ exponential complete Bell polynomials and $\{h_i\}$ cumulants of $-t \punt \alpha.$ Equation (\ref{(complBell)}) has been proved in \cite{Sole} by using Teugel martingales. 

For $Q_k(x,t),$ the Sheffer identity with respect to 
$t$ holds:
$$Q_k(x,t+s)=\sum_{j=0}^k \binom{k}{j} P_j(s) Q_{k-j}(x,t),$$
where $P_j(s)=Q_{j}(0,s)$ for all nonnegative integers $j.$ 
\subsection{Examples.}
The discretized version of a L\'evy process is a random walk $S_n = X_1 + X_2+ \cdots 
+ X_n,$ with $\{X_i\}$ independent and identically distributed r.v.'s.
For the symbolic representation of a L\'evy process we have dealt with, the symbolic counterpart of a random walk is the auxiliary umbra $n \punt \alpha.$ Indeed the 
infinite divisible property (\ref{(idp)}) is highlighted in the summation $\alpha^{\prime} + \alpha^{\prime \prime} + \cdots + \alpha^{\prime \prime \prime},$ encoded in the symbol $n \punt \alpha,$ with $\alpha^{\prime}, \alpha^{\prime \prime}, \ldots, \alpha^{\prime \prime \prime}$ uncorrelated and similar umbrae. Nevertheless not all r.v.'s having the symbolic representation $n \punt \alpha$ share the infinite divisible property. For example, the binomial r.v. has not 
the infinite divisible property \cite{sato}, nevertheless its symbolic representation is of type $n \punt \alpha$ where $\alpha \equiv \chi \punt p \punt \beta$ and $p \in (0,1).$ 
So the generality of the symbolic approach lies in the circumstance that if the parameter $n$ is replaced by $t,$ that is if the random walk is replaced by a L\'evy process, more general classes of polynomials can be recovered for which many of the properties here introduced still hold. The following tables resume the TSH representation for different families of classical polynomials, see \cite{Roman}. In particular Table 3 gives the umbra corresponding to the r.v. $X_i$ of $S_n$ in the first column,  its umbral counterpart in the second column and the associated TSH polynomial in the third column. In Table 4, the TSH polynomials given in Table 3 are traced back to special families of polynomials. In particular, with the polynomials ${\mathcal P}_k(x,t)$ we refer to 
$${\mathcal P}_k(x,t)=\sum_{j=1}^k Q_j(x,t) B_{k,j}(m_1, m_2, \ldots, m_{k-j+1})$$
for suitable $\{Q_j(x,t)\}$ and $\{m_i\}.$  
\begin{table}[ht]
\centering
\renewcommand{\arraystretch}{1.2}
\begin{tabular}{|l||c|c|} \hline
$X_i$        & Umbral counterpart    & Corresponding TSH polynomial  \\ \hline
Uniform $[0,1]$    & $-1 \punt \iota$ &  $E[(x + n \punt \iota)^k]$ \\ \hline 
Bernoulli $p=\frac{1}{2}$    & $\frac{1}{2}(-1 \punt \epsilon + u)$ &  $E[(x + n \punt \left[ \frac{1}{2}(-1 \punt u + \epsilon)\right])^k]$ \\ \hline 
Bernoulli $p \in (0,1)$    & $\chi \punt p \punt \beta$ &  $E[(x - n \punt \chi \punt p \punt \beta)^k]$ \\ \hline 
Sum of $a \in {\mathbb N}$ & $a \punt (-1 \punt \iota)$ & $E[(x + (a\,n) \punt \iota)^k]$  \\
uniform r.v.'s on $[0,1]$ &  &  \\ \hline
\end{tabular}
\caption{TSH polynomials associated to special random walks}
\label{table3}
\end{table}
\begin{table}[ht]
\centering
\renewcommand{\arraystretch}{1.2}
\begin{tabular}{|l||l|l|} \hline
$X_i$   & Special families & Connection with  \\ 
        & of polynomials & TSH polynomials \\ \hline
Uniform $[0,1]$   & Bernoulli $B_k(x,n)$        & $B_k(x,t)=Q_k(x,t)$ \\ \hline
Bernoulli $p=1/2$ & Euler ${\mathcal E}_k(x,n)$ & ${\mathcal E}_k(x,n)=Q_k(x,t)$ \\ \hline
Bernoulli $p \in (0,1)$    & Krawtchouk ${\mathcal K}_k(x,p,n)$ & $(n)_k 
{\mathcal K}_k(x,p,n)= {\mathcal P}_k(x,t)$ \\ 
                           &                                    & 
                           $m_i=E[((-1 \punt \chi \punt p \punt \beta)^{<-1>})^i]$ \\ \hline
Sum of $a \in {\mathbb N}$ & pseudo-Narumi $N_k(x,an)$ & $k! \, N_k(x,an) = {\mathcal P}_k(x,t)$ \\ 
uniform r.v.'s on $[0,1]$ &                                    & 
                           $m_i=E[(u^{<-1>})^i]$ \\ \hline

\end{tabular}
\caption{Connection between special families of polynomials and TSH polynomials}
\label{table4}
\end{table}
\par
Next tables 5 and 6 give TSH polynomials for some special L\'evy processes. 
\begin{table}[ht]
\centering
\renewcommand{\arraystretch}{1.2}
\begin{tabular}{|l||c|c|} \hline
L\'evy process           & Umbral representation    & TSH polynomial $Q_k(x,t)$  \\ \hline
Brownian motion                &                                  &   \\
with variance $s^2$            & $t \punt \beta \punt (s \varsigma)$ &  $E[(x - t \punt \beta \punt (s \varsigma))^k]$ \\ \hline
Poisson process                &                                  &  \\ 
with parameter $\lambda$       & $(t \lambda) \punt \beta$        &  $E[(x - (t \lambda) \punt \beta)^k]$ \\ \hline
Gamma process                  &                                  &   \\ 
with scale parameter $1$       & $t \punt \bar{u}$                & $E[(x - t \punt \bar{u})^k]$ \\
and shape parameter $1$        &                                  & \\ \hline
Gamma process                  &                                  & \\ 
with scale parameter $\lambda$ & $(t \lambda) \punt \bar{u}$      & $E[(x - (t \lambda) \punt \bar{u})^k]$ \\
and shape parameter $1$        &                                  & \\ \hline
Pascal process                 &                                  & \\ 
with parameter $d=p/q$         & $t \punt \bar{u} \punt d \punt \beta$  & $E[(x - t \punt \bar{u} \punt d \punt \beta)^k]$ \\
and $p+q=1$                    &                                  & \\ \hline
\end{tabular}
\caption{TSH polynomials associated to special L\'evy processes}
\label{table5}
\end{table}
\begin{table}[ht]
\centering
\renewcommand{\arraystretch}{1.2}
\begin{tabular}{|l||l|l|} \hline
L\'evy process   & Special family & Connection with  \\ 
                 & of polynomials & TSH polynomials \\ \hline
Brownian motion           &                          &  \\ 
with variance $s^2$       & Hermite $H_k^{(s^2)}(x)$ & $H_k^{(s^2)}(x)=Q_k(x,t)$ \\ \hline  
Poisson process           &                          & $\tilde{C}_k(x, \lambda t)=$ \\ 
with parameter $\lambda$  & Poisson-Charlier $\tilde{C}_k(x, \lambda t)$ & $=\sum_{j=1}^k s(k,j) Q_k(x,t)$ \\ 
                          &                                 & with $s(k,j)$ Stirling  \\ 
                          &                                 & numbers of first kind \\ \hline
Gamma process             &                                     & \\ 
with scale parameter $1$  & Laguerre ${\mathcal L}_k^{t-k}(x)$  & $k! (-1)^k  {\mathcal L}_k^{t-k}(x)= Q_k(x,t)$ \\
and shape parameter $1$   & & \\ \hline
Gamma process             & &   \\                      
with scale parameter $\lambda$ & actuarial $g_k(x,t)$    &  $g_k(x,t)= {\mathcal P}_k(x,t)$ \\ 
and shape parameter $1$        &                                  & $m_i=E[((\chi \punt (- \chi))^{<-1>})^i]$ \\ \hline

Pascal process                 &  &  $(-1)^k (t)_k M_k(x,t,p)=$  \\ 
with parameter $d=p/q$         & Meixner polynomials                                   & $= {\mathcal P}_k(x,t)$ \\
and $p+q=1$                    & of first kind $M_k(x,t,p)$                                  & $m_i=E[(\chi \punt (-1 \punt \chi + \chi / p))^i]$ \\ \hline
\end{tabular}
\caption{Special families of polynomials and TSH polynomials}
\label{table6}
\end{table}
\subsection{Orthogonality of TSH polynomials.}
A special class of TSH polynomials is the one including the L\'evy-Sheffer polynomials, whose applications within orthogonal polynomials are given in   
\cite{Schoutens}. A sequence of polynomials $\{V_k(x,t)\}_{t \geq 0}$ \cite{ST} is a L\'evy-Sheffer system if its generating function is such that 
\begin{equation}
\sum_{k \geq 0} V_k(x,t)\frac{z^k}{k!} = \left(g(z)\right)^t \exp\{x u(z)\},
\label{genfunlevshef}
\end{equation}
where
$g(z)$ and $u(z)$ are analytic functions in a neighborhood of $z = 0,$ $u(0) = 0,$ $g(0) = 1,$ $u'(0) \neq 0$ 
and $1/g(\tau(z))$ is an infinitely divisible moment generating function, with
$\tau(z)$ such that $\tau(u(z)) = z.$ Assume $\alpha$ an umbra such that $f(\alpha,z)=g(z)$ and $\gamma$ an umbra such that $f(\gamma,z)=1+u(z).$ 
From (\ref{genfunlevshef}), the L\'evy-Sheffer polynomials are moments of $x \punt \beta \punt \gamma + t \punt \alpha:$
\begin{equation}
V_k(x,t)=E[(x \punt \beta \punt \gamma + t \punt \alpha)^k].
\label{(VKLS)}
\end{equation}
\begin{thm} \label{aaa}
The TSH polynomials $Q_k(x,t)$ are special L\'evy-Sheffer polynomials.
\end{thm}
The proof of Theorem \ref{aaa} is straightforward by choosing in (\ref{(VKLS)}) as umbra $\alpha$ its inverse $-1 \punt \alpha$ and as umbra $\gamma$
the singleton umbra $\chi.$ All the L\'evy-Sheffer polynomials possess the TSH property. Indeed the following theorem has been proved in \cite{elviraannali}.
\begin{thm}
The L\'evy-Sheffer polynomials $\{V_k(x,t)\}_{t \geq 0}$ are TSH with respect to
L\'evy processes umbrally represented by   $\{-t \punt \alpha \punt \beta \punt \gamma^{{\scriptscriptstyle <-1>}}\}_{t \geq 0}.$
\end{thm}
In particular, one has \cite{elviraannali}
\begin{equation}
V_k(x,t) = \sum_{i=0}^k E[(x + t \punt \beta \punt \kappa_{({\alpha \punt \beta \punt \gamma^{{\scriptscriptstyle <-1>}}})})^i] B_{k,i}(g_1, \ldots, g_{k-i+1}),
\label{(LScomb)}
\end{equation}
 where $g_i = E[\gamma^i],$ for all nonnegative $i,$ and $\kappa_{({\alpha \punt \beta \punt \gamma^{{\scriptscriptstyle <-1>}}})}$ is the cumulant umbra of $\alpha \punt \beta \punt \gamma^{{\scriptscriptstyle <-1>}},$ with $\gamma^{{\scriptscriptstyle <-1>}}$ the compositional inverse of the umbra $\gamma.$ When the umbra $\alpha$ is replaced by its inverse and the umbra $\gamma$ by the singleton umbra, since $\chi^{{\scriptscriptstyle <-1>}} \equiv \chi,$ the only contribution in the summation (\ref{(LScomb)}) is given by $i=k.$ So again equation (\ref{(LScomb)}) reduces to $Q_k(x,t)=E[(x - t \punt \alpha)^k]$ since $\kappa_{({\alpha \punt \beta \punt \gamma^{{\scriptscriptstyle <-1>}}})} \equiv \chi \punt -1 \punt \alpha.$

Within L\'evy-Sheffer polynomials, the L\'evy-Meixner polynomials are those orthogonal with respect to the L\'evy processes $-t \punt \alpha \punt \beta \punt \gamma^{{\scriptscriptstyle <-1>}},$ 
due to their TSH property. The orthogonal property is 
$$E \left[ V_n(-t \punt \alpha \punt \beta \punt \gamma^{{\scriptscriptstyle <-1>}},t) V_m(-t \punt \alpha \punt \beta \punt \gamma^{{\scriptscriptstyle <-1>}},t)\right] = c_m \delta_{n,m}.$$
According to \cite{Schoutens}, all the polynomials in Table 6 are orthogonal. Their measure of orthogonality corresponds to the L\'evy process $-t \punt \alpha$ since 
$\chi^{{\scriptscriptstyle <-1>}} \equiv \chi$ and $-t \punt \alpha \punt \beta \punt \gamma^{{\scriptscriptstyle <-1>}} \equiv -t \punt \alpha.$ 
\subsection{Kailath-Segall polynomials.}
Equivalence (\ref{(TSHbell)}) gives the connection between TSH polynomials and Kailath-Segall polynomials \cite{KS}, which is a different class of polynomials strictly related to L\'evy processes. Indeed, both have a representation in terms of partition umbra of a suitable polynomial umbra. Overlaps are removed  by suitably choosing the indeterminates. 

The $n$-th Kailath-Segall polynomial $P_n(x_1, \ldots, x_n)$ is a multivariable polynomial such that when the indeterminates are replaced by the sequence $X_t^{(1)}, \ldots,$ $ X_t^{(n)}$ of variations of a L\'evy process $X_t$ 
$$X_t^{(1)} = X_t, \,\,\, X_t^{(2)} = [X,X]_t, \,\,\, X_t^{(n)} = \sum_{s \geq t}(\Delta X_s)^n
\;\; n \geq 3,$$
its iterated integrals are recovered
$$P_t^{(0)} = 1, \,\,\, P_t^{(1)}=X_t, \,\,\, P_t^{(n)}=\int_0^t P_{s-}^{(n-1)} {\rm d}X_s, \;\; n \geq 2,$$
that is $P_t^{(n)} = P_n\left(X_t^{(1)}, \ldots, X_t^{(n)}\right).$ The following recursion formula is known as Kailath-Segall formula  
\begin{equation}
P_t^{(n)} = \frac{1}{n} \left( P_t^{(n-1)} X_t^{(1)} - P_t^{(n-2)} X_t^{(2)} + \cdots +
(-1)^{n+1} P_t^{(0)} X_t^{(n)} \right).
\label{(KS)}
\end{equation}
When $X_t^{(1)}, \ldots, X_t^{(n)}$ are replaced by the power sums $S_1, \ldots, S_n$ in the 
indeterminates $x_1, \ldots, x_k,$ according to formula (1.2) in \cite{Taqqu} and Theorem 3.1 in \cite{Bernoulli}, the corresponding polynomials 
$P_n(S_1, \ldots, S_n)$ are such that 
\begin{equation}
n! P_n(S_1, \ldots, S_n) = E[ (\beta \punt [(\chi \punt \chi) \sigma])^n],
\label{(KSPS)}
\end{equation}
where $\sigma$ is the power sum umbra representing $\{S_j\}$ and the $\chi$-cumulant umbra $\chi \punt \chi$ represents the sequence $\{(-1)^{i-1} (i-1)!\}.$

In order to recognize special TSH polynomials within the family $\{P_n\},$ two steps are necessary: 
\begin{quote} 
{\it i)} Kailath-Segall polynomials need to be umbrally represented when the power sums $\{S_j\}$
are replaced by the indeterminates $\{x_i\};$ \\
{\it ii)} the indeterminates $\{x_i\}$ need to be replaced by suitable terms involving $x$ and $t.$
\end{quote}
For the first step, we will use equation (\ref{(KSPS)}). Assume $p$ an umbra representing the sequence $\{E[(\chi \punt x_i \punt \beta)^i]\}.$ Observe that $E[(\chi \punt x_i \punt \beta)^i]=x_i$ for all nonnegative $i.$ Then from (\ref{(KSPS)}) one has $n!  P_n( x_1, \ldots,  x_n) = E[(\beta \punt [(\chi \punt \chi) p \,])^n]$ so that the generating function of $P_n$ is 
$$f\left(\beta \punt [(\chi \punt \chi) p \,],z \right) = \exp \left( \sum_{n \geq 1} \frac{(-1)^{n+1}}{n} z^n x_n \right),$$
see also \cite{Yablonski}. The strength of this symbolic representation essentially relies on the properties of the partition umbra $\beta \punt [(\chi \punt \chi) p \,]$ reproducing those of Bell polynomials. For example, the following
property of Kailath-Segall polynomials 
\begin{equation}
P_n(a x_1, a^2 x_2, \ldots, a^n x_n) = a^n P_n(x_1, x_2, \ldots, x_n), \quad a \in \Real
\label{(TSHBELL)}
\end{equation}
is proved by observing that $\beta \punt [a (\chi \punt \chi) p \,] \equiv a (\beta \punt [(\chi \punt \chi) p \,]).$
For the next step, we need to characterize the indeterminates $x_1, x_2, \ldots$ such that
\begin{equation}
\kappa_{(x \punt u)} \, \dot{+} \, \kappa_{(- t \punt \alpha)} \equiv  (\chi \punt \chi) p \Rightarrow E[(\kappa_{(x \punt u)})^n] + E[(\kappa_{(- t \punt \alpha)})^n] = (-1)^{n-1}(n-1)! \, x_n.
\label{(TSHKS)}
\end{equation}
In the following we show some examples of how to perform this selection. These results extend the connections between TSH polynomials and Kailath-Segall polynomials analyzed in \cite{Sole}.

{\bf Generalized Hermite polynomials:} \\ \indent 
Since  $E[(\kappa_{(x \punt u)})^i]= x \, \delta_{i,1}$ and 
$E[(\kappa_{(- t \punt \beta \punt (s \varsigma))})^i] = s^2 \, t \, \delta_{i,2},$
then 
$$k! P_k(x, s^2 t, 0, \ldots, 0) = H_k^{(t)}(x)$$
where $\sum_{k \geq 0} H_k^{(t)}(x) z^k/k! = \exp\{xz - t z^2/ 2\}.$

{\bf Poisson-Charlier polynomials:} \\ \indent 
Poisson-Charlier polynomials $\{\widetilde{C}_k(x,t)\},$ with generating function 
$$\sum_{k\geq 0} \widetilde{C}_k(x,t) \frac{z^k}{k!} = e^{-tz}(1+z)^x,$$
are umbrally represented by  
\begin{equation}
\widetilde{C}_k(x,\lambda t) = E[(x \punt \chi  - t \punt \lambda \punt u)^k],
\label{(PCTSH)}
\end{equation}
see \cite{elviraannali}. Nevertheless (\ref{(PCTSH)}) differs from the result of Theorem \ref{UTSH2}, the TSH property holds since $\{\widetilde{C}_k(x,\lambda t)\}$ are a linear combination of special $Q_k(x,t).$ Moreover representation (\ref{(PCTSH)}) allows us the connection
with Kailath-Segall polynomials, when the indeterminates $\{x_i\}$ are chosen such that $E[(\kappa_{(x \punt \chi)})^i] + E[(\kappa_{(- t \punt \lambda \punt u)})^i] = (-1)^{i-1}(i-1)! \, x_i.$ 
Since 
$$E[(\kappa_{(x \punt \chi)})^i] + E[(\kappa_{(- t \punt \lambda \punt u)})^i] = \left\{ \begin{array}{lc} 
x - t \lambda, & i=1, \\
(-1)^{i-1} (i-1)! \, x^i, & i \geq 2,
\end{array} \right.$$
then $k! P_k(x - t \lambda,x,x, \ldots)=\widetilde{C}_k(x,\lambda t).$ 

{\bf Laguerre polynomials:} \\ \indent Laguerre polynomials $\{\mathcal{L}_k^{t-k}(x)\}$ are TSH polynomials 
such that 
$$k! (-1)^k  {\mathcal L}_k^{t-k}(x) = E[(x - t \punt \bar{u})^k].$$
They can be traced back to
Kailath-Segall polynomials if the indeterminates $\{x_i\}$ are characterized by
$E[(\kappa_{(x \punt u)})^i] + E[(\kappa_{(- t \punt \bar{u})})^i] = (-1)^{i-1}(i-1)! \, x_i.$ 
Since $f(\kappa_{(- t \punt \bar{u})},z) = 1 + t \, \log (1 -z)$ then $E[(\kappa_{(- t \punt \bar{u})})^i] = - t (i-1)!.$
So 
$$P_k(x-t, t , - t, t, \ldots) = (-1)^k {\mathcal L}_k^{t-k}(x)$$
and from (\ref{(TSHBELL)}) we have
$P_k(t-x, t , t, \ldots) = {\mathcal L}_k^{t-k}(x).$  

{\bf Actuarial polynomials:} \\ \indent
The actuarial polynomials $g_k(x,t)$ are a linear combination of suitable TSH polynomials
$Q_k(x,t)$ (see Table 6) but they are moments of the polynomial umbra $\lambda t - x \punt \beta,$ that is $g_k(x,t) = E[( \lambda t - x \punt \beta)^k].$ In order to characterize the
connection with Kailath-Segall polynomials, the indeterminates $\{x_i\}$ need to be characterized by
$E[(\kappa_{(\lambda t \punt u)})^i] + E[(\kappa_{(- x \punt \bar{u})})^i] = (-1)^{i-1}(i-1)! \, x_i.$
As before $E[(\kappa_{(\lambda t \punt u)})^i] = \lambda \, t \, \delta_{i,1},$ instead $E[(\kappa_{(- x \punt \bar{u})})^i]=-x (i-1)!$ for all nonnegative integers $i$ as in the previous example. Therefore one has $k! (-1)^k P_k(x - \lambda t, x , x, \ldots) = g_k(x,t) .$ 

{\bf Meixner polynomials of first kind:} \\ \indent
Meixner polynomials of first kind $\{M_k(x,t,p)\}$ \cite{Schoutens} 
are a linear combination of suitable TSH polynomials $Q_k(x,t)$ (see Table 6) but they are moments of the following polynomial umbra
$$
(-1)^k (t)_k M_k(x,t,p) = E\left\{\left[x \punt \left(-1 \punt \chi + \frac{\chi}{p}\right) 
- t \punt \chi\right]^k\right\},
$$
which allows us to find the connection with Kailath-Segall polynomials. Indeed for all nonnegative integers $i$ we have
$$E\left[\left\{\kappa_{x \punt \left(-1 \punt \chi + \frac{\chi}{p}\right)}\right\}^i\right]= (-1)^{i-1} \, (i-1)! \, x \, \left( \frac{1}{p^i} - 1 \right)$$
and
$$E\left[\left\{\kappa_{(\chi \punt -t \punt \chi)} \right\}^i\right] = (-1)^{i-1} (i-1)! \, t.$$
Then Kailath-Segall polynomials give the Meixner polynomials $(-1)^k (t)_k M_k(x,t,p)$ by choosing
$$x_i = \left[ \left( \frac{1}{p^i} - 1 \right) x - t \right]$$
for $i=1,2,\ldots.$ 
\section{Symbolic multivariate L\'evy processes.}
In the multivariate case, the main device of the symbolic method here proposed relies on the employment of multi-indices of length $d$. Sequences like $\{g_{i_1, i_2, \ldots, i_d}\}$ are replaced with a product of powers $\mu_{1}^{i_1} \mu_{2}^{i_2} \cdots \mu_{d}^{i_d},$ where $(\mu_1, \mu_2, \ldots, \mu_d)$ are umbral monomials and $(i_1, i_2, \ldots, i_d)$ are nonnegative integers. Since the umbral monomials could not have disjoint support, then the evaluation $E$ does not necessarily factorizes on the product $\mu_{1}^{i_1} \mu_{2}^{i_2} \cdots \mu_{d}^{i_d},$ that is
\begin{equation}
E[\mu_1^{i_1} \mu_2^{i_2} \cdots \mu_d^{i_d}]=E[\mubf^{\ibs}] = g_{\ibs}
\label{(multiindex)}
\end{equation}
where $\ibs=(i_1, i_2, \ldots, i_d)$ and $\mubf=(\mu_1, \mu_2, \ldots, \mu_d).$
We assume $g_{\boldsymbol 0}=1$ with ${\boldsymbol 0}=(0,0,\ldots,0).$ Then
$g_{\ibs}$ is called the multivariate moment of $\mubf.$ Table 7 shows some special $d$-tuples we will use later. 

\begin{table}[ht]
\centering
\renewcommand{\arraystretch}{1.1}
\begin{tabular}{|c|c|l|} \hline
\emph{ } & \emph{$d$-tuple} & \emph{Generating functions} \\ \hline
Multivariate Unity  $\ubs$ &  $(u, \ldots, u^{\prime})$ & $f(\ubs,\zbs) = e^{z_1 + \cdots + z_d}.$ \\ \hline
Multivariate Gaussian $\deltabf$ & $(\varsigma, \ldots, \varsigma^{\prime})$ & $f(\deltabf,\zbs)=1 + \frac{1}{2} \zbs \zbs^T.$ \\ \hline
Multivariate Bernoulli $\iotabf$ & $(\iota, \ldots, \iota)$ & $f(\iotabf,\zbs)=\displaystyle{\frac{z_1 + \cdots + z_d}{e^{z_1 + \cdots + z_d} - 1}}.$ \\ \hline
Multivariate Euler $\boeta$ & $(\eta, \ldots, \eta)$ & $f(\boeta,\zbs)=\displaystyle{\frac{2 e^{(z_1 + \cdots + z_d)}}{e^{2(z_1 + \cdots + z_d)} + 1}}.$ \\ \hline
\end{tabular}
\caption{Generating functions of special $d$-tuples of umbral monomials}
\label{table7}
\end{table}

The notions of similarity and uncorrelation are updated as follows. 
Two $d$-tuples $\mubf$ and $\nubs$ of umbral monomials are said to be similar
if they represent the same sequence of multivariate moments. They are said to
be uncorrelated if $E[\mubf^{{\ibs}_1} \nubs^{{\ibs}_2}] = E[\mubf^{{\ibs}_1}] 
E[\nubs^{{\ibs}_2}].$ 

Multivariate L\'evy processes are represented by $d$-tuples of umbral monomials.  
\begin{defn} \label{def_levy_multi}
A stochastic process $\{\X_t\}_{t \geq 0}$ on $\Real^d$ is a \emph{multivariate L\'evy process} if
\begin{description}
	\item[{\rm (i)}] $\X_0 = \boldsymbol{0}$ a.s.
	
	\item[{\rm (ii)}] For all $n \geq 1$ and for all $\,0 \leq t_1 \leq t_2 \leq \ldots \leq t_n < \infty,$ the r.v.'s $\X_{t_2} - \X_{t_1}, \X_{t_3} - \X_{t_2}, \ldots$ are independent.
	
	\item[{\rm (iii)}] For all $s \leq t,$ $\X_{t + s} - \X_s \stackrel{d}{=} \X_t.$
	
	\item[{\rm (iv)}] For all $\varepsilon > 0,$ $\lim_{h \rightarrow 0} P(|\X_{t + h} - \X_t| > \varepsilon) = 0.$
	
	\item[{\rm (v)}] $t \mapsto \X_t (\omega)$ are right-continuous with left limits, for all $\omega \in \varOmega,$ with $\varOmega$ the underlying sample space.   
\end{description}
\end{defn}
As in the univariate case, the moment generating function of a multivariate L\'evy process is $\varphi_{\scriptscriptstyle{\X_1}}(\zbs) = E\left[e^{\zbs \X_1^{\trasp}}\right],$ with $\zbs \in \Real^d.$ Paralleling the univariate case, the generating function of a $d$-tuple $\mubf$ is 
$$f(\mubf, \zbs) = 1 + \sum_{k \geq 1} \sum_{\substack{\ibs \in \mathbb{N}_0^d \\ |\ibs| = k}} g_{\ibs} \frac{\zbs^{\ibs}}{\ibs!}.$$
Choose the $d$-tuple $\mubf$ such that $f(\mubf,\zbs)=\varphi_{\X_1}(\zbs),$ that is
$E[\mubf^{\ibs}]=E[\X_1^{\ibs}]$ for all $\ibs \in \mathbb{N}_0^d.$ 
The auxiliary umbra $n \punt \mubf$ denotes the sum of $n$ uncorrelated $d$-tuples of umbral monomials similar to $\mubf.$
Its multivariate moment is \cite{dibruno}
\begin{equation}
E[(n \punt \mubf)^{\ibs}] = \sum_{\lambs \vdash \ibs} \frac{\ibs!}{\mathfrak{m}(\lambs) \lambs!} \, (n)_{l(\lambs)} \, 
E[\mubf_{\lambs}], \label{(auxmult)}
\end{equation}
where $E[\mubf_{\scriptscriptstyle{\lambs}}] = g_{\lambs_1}^{r_1} 
g_{\lambs_2}^{r_2} \ldots$ and $\lambs$ is a partition\footnote{A partition $\lambs$ of a multi-index $\ibs,$ in symbols $\lambs \vdash \ibs,$ is a 
matrix $\lambs = (\lambda_{ij})$ of nonnegative integers and with no zero columns in lexicographic order
$\prec$ such that $\lambda_{r_1} + \lambda_{r_2} + \cdots + \lambda_{r_k} = i_r$ for $r = 1, 2, \ldots , d.$ 
The number of columns of $\lambs$ is denoted by $l(\lambs)$. The notation $\lambs
= (\lambs_{1}^{r_1} , \lambs_{2}^{r_2}, \ldots)$ represents the
matrix $\lambs$ with $r_1$ columns equal to $\lambs_{1},$ $r_2$ columns equal to $\lambs_{2}$ and so on, 
where $\lambs_{1} \prec  \lambs_{2} \prec \ldots.$ We set $\mathfrak{m}(\lambs) = (r_1, r_2, \ldots),$ $\mathfrak{m}(\lambs)! = r_1! r_2! \cdots$ and $\lambs! = \lambs_1! \lambs_2! \cdots.$} 
of the multi-index $\ibs$ of length $l(\lambs).$ By replacing the nonnegative integer $n$ with the real parameter $t$ in (\ref{(auxmult)}) the resulting auxiliary
umbra $t \punt \mubf$ is the symbolic representation of the multivariate L\'evy process $ \X_t.$ 

As in the univariate case, since $\mubf \equiv \beta \punt \kappa_{\mubf}$ with $\kappa_{\mubf}$ the $\mubf$-cumulant umbra \cite{Bernoulli},
a different representation for a multivariate L\'evy process is $t \punt \beta \punt \kappa_{\mubf}.$ The cumulant $d$-tuple
could be further specified by using the multivariate L\'evy-Khintchine formula \cite{sato}. 
\begin{thm} \label{LK_thm_multi}
$\X = \{\X_t\}_{t \geq 0}$ is a L\'evy process if and only if there exists $\mbs_1 \in \Real^d,$ a symmetric, positive defined $d \times d$ matrix $\Sigma > 0$ and a measure $\nu$ on $\Real^d$ with
$$\nu(\{0\}) = 0 \mbox{ and } \int_{\Real}(|\xbs|^2 \wedge 1) \nu(d \xbs) < \infty$$
such that
\begin{equation}
\varphi_{\scriptscriptstyle{\X}}(\zbs) = \exp\left\{t\left[\frac{1}{2} \zbs \Sigma \zbs^{T} + \mbs_1\zbs^{T} + \int_{\Real^d}
(e^{\xbs\zbs^{T}} - 1 - \xbs\zbs^{T} {\boldsymbol 1}_{\{|\xbs|\leq 1\}}(\xbs))\, \nu({\rm d}\xbs)\right]\right\}.
\label{LK_multi}
\end{equation}
The representation of $\varphi_{\scriptscriptstyle \X}(\zbs)$ in (\ref{LK_multi}) by $\boldsymbol{m_1},$ $\Sigma$ and $\nu$ 
is unique.
\end{thm}
Set $\mbs_2 \zbs^{\trasp} = \int_{\Real^d}\zbs\xbs ^{\trasp} {\boldsymbol 1}_{\{|\xbs| > 1\}}(\xbs)\, \nu(d\xbs)$ and $\mbs = \mbs_1 
+ \mbs_2,$ then 
\begin{equation*}
\varphi_{\scriptscriptstyle{\X}}(\zbs) = \exp\bigg\{t\bigg[\frac{1}{2}\zbs\Sigma \zbs^{\trasp} + \mbs \zbs^{\trasp} + \int_{\Real^d}(e^{\zbs\xbs^{\trasp}} - 1 -\zbs\xbs^{\trasp}) \, \nu(d\xbs)\bigg]\bigg\},
\end{equation*}
that is,
\begin{equation}
\varphi_{\scriptscriptstyle{\X}}(\zbs) = \exp\bigg\{t\bigg[\frac{1}{2}\zbs\Sigma \zbs^{\trasp} + \mbs \zbs^{\trasp}\bigg]\bigg\}  \exp\bigg\{t\bigg[\int_{\Real^d}(e^{\zbs\xbs^{\trasp}} - 1 -\zbs\xbs^{\trasp}) \, \nu(d\xbs)\bigg]\bigg\}.
\label{LK2_multi}
\end{equation}
\begin{thm} \label{fund}
Every L\'evy process $\{\X_t\}_{t \geq 0}$ on $\Real^d$ is umbrally represented by the family of auxiliary umbrae
\begin{equation}
 \{t \punt \beta \punt (\chi \punt \mbs \dot{+} \deltabf C^{T} \dot{+} \gabs)\}_{t \geq 0},\label{levymulti}
\end{equation}
where $\beta$ is the Bell umbra, $\mbs \in \Real^d,$ $\deltabf$ is the multivariate umbral counterpart of a standard gaussian r.v., $C$ is the square root of the covariance matrix $\Sigma$ and $\gabs$ is the multivariate umbra associated 
to the L\'evy measure.
\end{thm}
Every auxiliary umbra $t \punt \beta \punt \kappa_{\mubf}$ is the symbolic version of a multivariate compound Poisson r.v. 
of parameter $t,$  that is a random sum $S_N = \Y_1 + \cdots + \Y_N$ of independent and identically distributed random vectors 
$\{ \Y_i\},$ whose index $N$ is a Poisson r.v. of parameter $t.$ Then the same holds for the L\'evy process. The $d$-tuple 
$(\deltabf C^{T}  \dot{+} \chi \punt \mbs \dot{+} \gabs)$ umbrally  represents any of the random vectors $\{ \Y_i\}.$ Observe that 
$\chi \punt \mbs$ has not a probabilistic counterpart. If $\mbs$ is not equal to the zero vector, this parallels the well-known 
difficulty to interpret the L\'evy  measure as a probability measure.
\subsection{Multivariate TSH polynomials.}
The conditional evaluation with respect to an umbral $d$-tuple $\mubf$ has been introduced in \cite{DiNardoOliva3}. 
Assume ${\mathcal{X}} = \{\mu_1, \mu_2, \ldots, \mu_d\}.$ The conditional evaluation with respect to the umbral $d$-tuple $\mubf$
is the linear operator 
$$E(\;\cdot \; \vline \,\, \mubf): \, \Real[x_1, \ldots, x_d][\A] \; \longrightarrow \; \Real[\mathcal{X}]$$ 
such that $E(1 \,\, \vline \,\, \mubf) = 1$ and 
\begin{equation}
E(x_1^{l_1} \, x_2^{l_2} \ \cdots \, x_d^{l_d} \, \mubf^{\ibs} \, \nubs^{\jbs} \, \gabs^{\kbs} \cdots \, \,  \vline \,\, \mubf) 
= x_1^{l_1} \, x_2^{l_2} \, \cdots \, x_d^{l_d} \, \mubf^{\ibs} \, E[\nubs^{\jbs}] \, E[\gabs^{\kbs}] \cdots
\label{(condeval1)}
\end{equation}
for uncorrelated $d$-tuples $\mubf, \nubs, \gabs \ldots,$ multi-indices $\ibs, \jbs, \kbs \ldots \in \mathbb{N}_0^d$ and $\{l_i\}_{i=1}^d$ nonnegative integers.
Since $f[(n + m) \punt \mubf, \zbs] = f(\mubf, \zbs)^{n + m} = f(n \punt \mubf, \zbs) \, f(m \punt \mubf, \zbs),$ then 
$(n + m) \punt \mubf \equiv n \punt \mubf + m \punt \mubf^{\prime},$ with $\mubf$ and $\mubf^{\prime}$ uncorrelated $d$-tuples of umbral monomials. Then, for 
$E(\;\cdot \; \vline \,\, n \punt \mubf)$ we assume  
$E[\{(n + m) \punt \mubf\}^{\ibs} \,\, \vline \,\, n \punt \mubf] = E[\{n \punt \mubf + m \punt \mubf^{\prime}\}^{\ibs} \,\, 
\vline \,\, n \punt \mubf]$ for all nonnegative integers $n, m$ and for all $\ibs \in \mathbb{N}_0^d.$ If $n \ne m,$ then
\begin{equation}
E[\{(n + m) \punt \mubf\}^{\ibs} \,\, \vline \,\, n \punt \mubf] = E[\{n \punt \mubf + m \punt \mubf\}^{\ibs} \,\, 
\vline \,\, n \punt \mubf],
\label{(aaa1)}
\end{equation}
since $n \punt \mubf$ and $m \punt \mubf$ are uncorrelated auxiliary umbrae. We will use the same $d$-tuple $\mubf$ as in (\ref{(aaa1)}) when no misunderstanding occurs.  Thanks to equations (\ref{(condeval1)}) and (\ref{(aaa1)}), 
we have
\begin{equation}
E\left[\{(n + m) \punt \mubf\}^{\ibs} \,\, \vline \,\, n \punt \mubf\right] = \sum_{\kbs \leq \ibs} \binom{\ibs}{\kbs} (n \punt \mubf)^{\kbs} E[(m \punt \mubf)^{\ibs - \kbs}],
\label{ce_dot_multi}
\end{equation}
where $\kbs \leq \ibs \; \Leftrightarrow \; k_j \leq i_j \mbox{ for all } j = 1, \ldots, d$ and 
$\binom{\kbs}{\ibs} = \binom{k_1}{i_1} \cdots \binom{k_d}{i_d}.$
By analogy with (\ref{(aaa1)}) and (\ref{ce_dot_multi}), we have $t \punt \mubf \equiv s 
\punt \mubf + (t-s) \punt \mubf$ and for $t \geq 0$ and $s \leq t$
$$
E\left[(t \punt \mubf)^{\ibs} \,\, \vline \,\, s \punt \mubf \right] = \sum_{\kbs \leq \ibs} \binom{\ibs}{\kbs} (s \punt \mubf)^{\kbs} 
E[\{(t - s) \punt \mubf\}^{\ibs - \kbs}].
$$
\begin{thm}\label{UTSH2_multi}
For all $\ibs \in \mathbb{N}_0^d,$ the family of polynomials
\begin{equation}
Q_{\ibs}(\xbs,t) = E[(\xbs - t \punt \mubf)^{\ibs}] \in \Real[x_1, \ldots, x_d]
\label{(tshumbral_multi)}
\end{equation}
is TSH with respect to $\{t \punt \mubf\}_{t \geq 0}.$
\end{thm}
The auxiliary umbra $-t \punt \mubf$ denotes the inverse of $t \punt \mubf$ that is $- t \punt \mubf + t \punt \mubf \equiv \varbs,$ where $\varbs$ is the $d$-tuple such that $\varbs=(\epsilon_1, \epsilon_2, \ldots, \epsilon_d)$ with $\{\epsilon_i\}$ uncorrelated augmentation umbrae. Coefficients of $Q_{\ibs}(\xbs,t)$ in (\ref{(tshumbral_multi)}) are such that
$$Q_{\ibs}(\xbs,t) = \sum_{\kbs \leq \ibs} \binom{\ibs}{\kbs} \xbs^{\ibs - \kbs} E[(-t \punt \mubf)^{\kbs}]$$
so when $\xbs$ is replaced by $t \punt \mubf$ their overall evaluation is zero. Properties on the coefficients of $Q_{\ibs}(\xbs,t)$ can be found
in \cite{DiNardoOliva3}. Here we just recall a characterization of the coefficients of any multivariate TSH polynomial in terms
of those of $Q_{\ibs}(\xbs,t).$ 
\begin{thm} \label{comb_multi}
A polynomial 
\begin{equation}
P(\xbs,t) =\sum_{\kbs \leq \vbs} p_{\scriptscriptstyle\kbs}(t) \, \xbs^{\kbs}
\label{polynomial}
\end{equation}
is a TSH polynomial with respect to $\{t \punt \mubf\}_{t \geq 0}$ if and only if
\begin{equation}
p_{\scriptscriptstyle\kbs}(t) = \sum_{\kbs \leq \ibs \leq \vbs}  \binom{\ibs}{\kbs} \, p_{\scriptscriptstyle\kbs}(0) \, 
E[(-t \punt \mubf)^{\ibs-\kbs}], \quad \hbox{ for } \kbs \leq \vbs. 
\label{prop5_multi}
\end{equation}
\end{thm}

Table 8 and 9 give some examples of multivariate TSH polynomials
and their connection with multivariate L\'evy processes. The corresponding $d$-tuples are given in Table 7.
\begin{table}[ht]
\centering
\renewcommand{\arraystretch}{1.2}
\begin{tabular}{|l||c|c|} \hline
multivariate L\'evy process                 & Umbral representation                                & TSH polynomial $Q_{\ibs}(\xbs,t)$  \\ \hline
Brownian motion                             &                                                      &   \\
with covariance $\Sigma=C C^T$              & $t \punt \beta \punt (\deltabf C^T)$                 &  $E[(\xbs - t \punt \beta \punt (\deltabf C^T))^{\ibs}]$ \\ \hline
$\X_t$ with $\X_1 \stackrel{d}{=} (U, \ldots, U)$  &                                                      &  \\ 
and $U$ uniform r.v. on $[0,1]$             & $- t \punt \iotabf$                                  &  $E[(\xbs + t \punt \iotabf)^{\ibs}]$ \\ \hline
$\X_t$ with $\X_1 \stackrel{d}{=} (Y, \ldots, Y)$  &                                                      &   \\ 
and $Y$ Bernoulli r.v.                      & $\frac{1}{2} [t \punt (\ubs -  1 \punt \boeta)]$     &  $E\left\{ \left(\xbs + \frac{1}{2} [ t \punt (\boeta - \ubs)] \right)^{\ibs}\right\}$ \\
of parameter $1/2$                          &                                  & \\ \hline
\end{tabular}
\caption{TSH polynomials associated to special multivariate L\'evy processes}
\label{table8}
\end{table}
\begin{table}[ht]
\centering
\renewcommand{\arraystretch}{1.2}
\begin{tabular}{|l||l|l|} \hline
multivariate L\'evy process                          & Special family                             & Connection with  \\ 
                                                     & of polynomials                             & TSH polynomials \\ \hline
Brownian motion                                      &                                            &  \\ 
with covariance $\Sigma=C C^T$                       & Hermite $H_{\ibs}^{(t^2)}(\xbs,\Sigma)$    & $H_{\ibs}^{(t^2)}(\xbs,\Sigma)=Q_{\ibs}(\xbs,t)$ \\ \hline  
$\X_t$ with $\X_1 \stackrel{d}{=} (U, \ldots, U)$           &                                            &   \\ 
and $U$ uniform r.v. on $[0,1]$                      & Bernoulli $B_{\ibs}^{(t)}(\xbs)$           & $B_{\ibs}^{(t)}(\xbs) = Q_{\ibs}(\xbs,t)$ \\ \hline
$\X_t$ with $\X_1 \stackrel{d}{=} (Y, \ldots, Y)$           &                                            & \\ 
and $Y$ Bernoulli r.v.                               & Euler $\mathcal{E}_{\ibs}^{(t)}(\xbs)$     & $\mathcal{E}_{\ibs}^{(t)}(\xbs) = Q_{\ibs}(\xbs,t)$ \\
of parameter $1/2$                                   &                                            & \\ \hline
\end{tabular}
\caption{Special families of polynomials and TSH polynomials}
\label{table9}
\end{table}

Let us remark that Hermite polynomials  $H_{\ibs}^{(t^2)}(\xbs,\Sigma)$ in Table 9 are a generalization of the polynomials  $H_{\ibs}(\xbs)$ in \cite{Withers} 
whose moment representation is $H_{\ibs}(\xbs)=E[(\xbs \Sigma^{-1} + {\mathfrak i} \Ybs)^{\ibs}]$ 
with $E$ the expectation symbol, $\Ybs \simeq N({\bf 0}, \Sigma^{-1})$ and $\Sigma$ a covariance matrix of full rank $d$.

A generalization of L\'evy-Sheffer system to the multivariate case has been introduced in \cite{DiNardoOliva3}. A sequence of multivariate polynomials $\{V_{\kbs}(\xbs,t)\}_{t \geq 0}$ is a multivariate L\'evy-Sheffer system if 
$$
1 + \sum_{k \geq 1} \sum_{\substack{\vbs \in \mathbb{N}_0^d \\ |\vbs| = k}} V_{\kbs}(\xbs,t)\frac{\zbs^{\kbs}}{\kbs!} =  [g(\zbs)]^t 
\exp\{ (x_1 + \cdots + x_d) [h(\zbs) - 1] \},
$$
where $g(\zbs)$ and $h(\zbs)$ are analytic in a neighborhood of $\zbs =\boldsymbol{0}$ and 
$$\left. \frac{\partial}{\partial z_i} h(\zbs) \right|_{\zbs = \boldsymbol{0}} \neq 0 \quad \hbox{for 
$i=1,2,\ldots,d$}.$$ 
If $\mubf$ and $\nubs$ are $d$-tuples of umbral monomials such that $f(\mubf,\zbs)=g(\zbs)$ and $f(\nubs,\zbs)=h(\zbs)$ respectively, then 
\begin{equation}
V_{\kbs}(\xbs,t) = E[(t \punt \mubf +  (x_1 + \cdots + x_d) \punt \beta \punt \nubs)^{\kbs}].
\label{(LS)}
\end{equation}
The multivariate L\'evy-Sheffer polynomials  for the pair $\mubf$  and $\nubs$ are TSH polynomials with respect to a special symbolic multivariate L\'evy process
involving the multivariate compositional inverse of a $d$-tuple $\nubs.$  Assume  $\chibs_{(i)}$ the $d$-tuple with all components equal to the augmentation umbra and only the $i$-th one equal
to the singleton umbra, that is $\chibs_{(i)} = (\epsilon, \ldots, \chi, \ldots, \epsilon).$
The multivariate compositional inverse of $\nubs$ is the umbral $d$-tuple 
$\nubs^{{\scriptscriptstyle <-1>}}=({(\nubs^{{\scriptscriptstyle <-1>}})}_1, \ldots, (\nubs^{{\scriptscriptstyle <-1>}})_d)$ 
such that $(\nubs^{{\scriptscriptstyle <-1>}})_i \punt \beta \punt \nubs \equiv \chibs_{(i)}$ for $i = 1, \ldots, d.$ 

\begin{thm}
The multivariate L\'evy-Sheffer polynomials  for the pair $\mubf$  and $\nubs$ are TSH polynomials with respect to the symbolic multivariate L\'evy process 
$$\{t \punt (\mu_1 \punt \beta \punt \nubs^{{\scriptscriptstyle <-1>}}_1 + \cdots + \mu_d \punt \beta \punt \nubs^{{\scriptscriptstyle <-1>}}_d)\}_{t \geq 0}.$$
\end{thm}
\section{Conclusions and open problems.}
In this paper, the review of a symbolic treatment of TSH polynomials,
relied on the classical umbral calculus, is proposed. The main advantage of this symbolic presentation is the plainness of the overall setting which reduces to few fundamental statements,
but also the availability of efficient routines \cite{dinardooliva2009} for the implementation of formulae as (\ref{(auxmult)}), which is the key to manage
the polynomials $Q_{\kbs}(\xbs,t).$

The main result of this presentation is that any univariate (respectively multivariate) TSH polynomial has the form $Q_k(x,t)$ (respectively 
$Q_{\kbs}(\xbs,t)$) or can be expressed as a linear combination of the polynomials $Q_k(x,t)$ with coefficients given by (\ref{prop5_multi}).  Thanks to the umbral
representation of multivariate L\'evy-Sheffer systems, more families of umbral polynomials could be characterized, together with their orthogonality properties. 
This will be the object of future research  and investigation.  

In \cite{Barrieu}, Barrieu and Shoutens have related the infinitesimal generator of a Markov process to a more general class of linear operators possessing the TSH property, both ascribable to special families of martingales. A stochastic Taylor formula is produced which results to be a generalization of a TSH polynomial due to the presence of a remainder term series. A symbolic representation of this new TSH function could open the way to a new classification of the 
corresponding operators by which to recover the martingale property on L\'evy processes. Similarly, the extension to the more general class of Markov processes (a first attempt is given in \cite{Barrieu})
would move the employment of TSH functions beyond the field of applications strictly connected to the market portfolio. 
One step more consists in dealing with matrix-valued stochastic processes by replacing formal power series (\ref{(fps)}) with hypergeometric functions, as done in 
\cite{lawi}. This would allow us a symbolic representation also for zonal polynomials whose computational handling is still an open problem.

\end{document}